\documentclass[a4paper,12pt,reqno]{amsart}
\usepackage{amssymb}
\usepackage[dvips,]{graphicx}

\input epsf.tex 
\newdimen\xsize
\newdimen\oldbaselineskip
\newdimen\oldlineskiplimit
\def\restorelineskip{\baselineskip=\oldbaselineskip%
\lineskiplimit=\oldlineskiplimit}
\def\putm[#1][#2]#3{
\hbox{\vbox to 0pt{\parindent=0pt%
\vskip#2\xsize\hbox to0pt{\hskip#1\xsize $#3$\hss}\vss}}}%
\long\def\Line#1{\hbox to \hsize{#1}}
\def\putt[#1][#2]#3{
\vbox to 0pt{\noindent\hskip#1\xsize\lower#2\xsize%
\vtop{\restorelineskip#3}\vss}}

\makeatletter
\def\xbig[#1]#2{{\hbox{$\m@th\left#2\vbox to#1\xsize{}%
\right.\n@space$}}}
\makeatother
\def\xlar[#1]#2{%
\smash{\mathop{ \hbox to #1\xsize{\leftarrowfill}}\limits^{#2}}}
\def\xrar[#1]#2{%
\smash{\mathop{ \hbox to #1\xsize{\rightarrowfill}}\limits^{#2}}}
\def\xline[#1]{\hbox to #1\xsize{\leaders\hrule\hfill}}

\DeclareFontFamily{U}{rsf}{\skewchar\font'177}%
\DeclareFontShape{U}{rsf}{m}{n}{<-6>rsfs5<6-8>rsfs7<8->rsfs10}{}%
\DeclareFontShape{U}{rsf}{b}{n}{<-6>rsfs5<6-8>rsfs7<8->rsfs10}{}%
\DeclareMathAlphabet\RSFS{U}{rsf}{m}{n}
\SetMathAlphabet\RSFS{bold}{U}{rsf}{b}{n}
\def\sf#1{{\mathsf{#1}}}

\def\slsf{\slshape \sffamily }

\def\msmall#1{\mathchoice{\hbox{\small$\displaystyle {#1}$}}{#1}{#1}{#1}}

\def\cc{{\mathbb C}}
\def\dd{{\mathbb D}}

\def\rr{{\mathbb R}}

\def\lim{\mathop{\sf{lim}}}

\def\vol{\sf{Vol}}

\def\<{\langle}\let\la=\<
\def\>{\rangle}\let\ra=\>
  
\def\comp{\Subset}
\def\d{\partial}

\def\ddef{\mathrel{{=}\raise0.3pt\hbox{:}}}
\def\deff{\mathrel{\raise0.3pt\hbox{\rm:}{=}}}

\def\fraction#1/#2{\mathchoice{{\msmall{ #1\over#2}}}%
{{ #1\over #2 }}{{#1/#2}}{{#1/#2}}}

\def\le{\leqslant}

\def\emptyset{\varnothing}

\def\longpoints{\leaders\hbox to 0.5em{\hss.\hss}\hfill \hskip0pt}
\def\stateskip{\smallskip}
\def\state#1. {\stateskip\noindent{\bf#1. }} 
\def\statep#1. {\stateskip\noindent{\bf#1 }} 
\def\proof{\state Proof. \2}

\def\Chi{\raise 2pt\hbox{$\chi$}}

\def\Chi{\raise 2pt\hbox{$\chi$}}

\let\phI=\phi\let\phi=\varphi\let\varphi=\phI
\let\cal=\mathcal

\def\calf{{\cal F}}

\def\comp{\Subset}
\def\d{\partial}

\def\1{{1\mkern-5mu{\rom l}}}

\def\ge{\geqslant}

\def\fraction#1/#2{\mathchoice{{\msmall{ #1\over#2}}}%
{{ #1\over #2 }}{{#1/#2}}{{#1/#2}}}

\def\le{\leqslant}

\def\emptyset{\varnothing}
\newcommand{\2}{\thinspace}

\def\qed{\ \ \hfill\hbox to .1pt{}\hfill\hbox to .1pt{}\hfill $\square$\par}

\def\comment#1\endcomment{}

 \belowdisplayskip=\abovedisplayskip

\def\lineeqqno(#1){\hfill\llap{\vbox to 10pt%
{\vss\begin{align} \eqqno(#1)\end{align}\vss}}\vskip1pt}

\advance\headheight 1.2pt

\def\ShowwLLabel#1{}

\def\thechpt{\Roman{chpt}}

\def\newchapt[#1]#2{\newpage%
\refstepcounter{chpt}\setcounter{subsection}{0}%
\setcounter{thm}{0}\setcounter{defi}{0}%
\setcounter{rema}{0}\setcounter{exrc}{0}%
\renewcommand{\thesubsection}{\thechpt.\arabic{subsection}}%
\section*{\begin{center}\huge \bf Chapter \thechpt\\
#2 \end{center}}\label{#1}%
\ \smallskip%
\markboth{Chapter \thechpt}{#2}%
}
\def\newsect[#1]#2{\refstepcounter{section}\setcounter{equation}{0}%
\renewcommand{\thesubsection}{\arabic{section}.\arabic{subsection}}%
\section*{\arabic{section}.
#2}\vspace{-20pt}\label{#1}\vspace{20pt}%
\markboth{Section \arabic{section}}{#2}}

\def\newlect[#1]#2{\refstepcounter{section}%
\renewcommand{\thesubsection}{\arabic{section}.\arabic{subsection}}%
\section*{Lecture \arabic{section}\\
#2}\label{#1}%
\markboth{Lecture \arabic{section}}{#2}}

\def\newprg[#1]#2{\refstepcounter{subsection}%
\subsection*{{\thesubsection.\ #2}} \label{#1}%
}

\def\newappx[#1]#2{%
\refstepcounter{appx}\setcounter{section}{0}%
\renewcommand{\thesubsection}{A\arabic{appx}.\arabic{subsection}}%
\section*{Appendix \arabic{appx}\\ #2}
\label{#1}%
\markboth{Appendix A\arabic{appx}}{#2}
}

\newtheorem{thm}{Theorem}
   \def\newthm#1{\begin{thm}\label{#1}}

\newtheorem{nnthm}{Theorem.} 
   \def\newnnthm#1{\begin{nnthm} \label{#1}}
\newtheorem{lem}{Lemma}
   \def\newlemma#1{\begin{lem} \label{#1}}

\newtheorem{prop}{Proposition}
   \def\newprop#1{\begin{prop}\label{#1}}

\newtheorem{corol}{Corollary}
   \def\newcorol#1{\begin{corol} \label{#1}}

\newtheorem{defi}{Definition}
   \def\newdefi#1{\begin{defi} \label{#1}\rm }

\newtheorem{exmp}{Example}
   \def\newexmp#1{\begin{exmp} \label{#1}\rm }

\newtheorem{exrc}{Exercise}
   \def\newexrc#1{\begin{exrc} \label{#1}\rm }

\newtheorem{rema}{Remark}
   \def\newrema#1{\begin{rema} \label{#1}\rm }

\def\eqqno(#1){\label{(#1)}}
\def\eqqref(#1){(\ref{(#1)})}

\title{Equicontinuous families of meromorphic mappings with values in compact
complex surfaces}
\author {F. Neji}
\date{\today}
\address{
Universit\'e de Lille-1, UFR de Math\'ematiques, 59655 Villeneuve
d'Ascq, France} \email{Fethi.neji@math.univ-lille1.fr}
\begin{document}

\begin{abstract}
We prove that a family of meromorphic mappings from a bidisc to a
compact complex surface, which are equicontinuous in a neighborhood of the boundary of the bidisc, has the volumes of its graphs locally uniformly bounded.
\end{abstract}

\maketitle

\newsect[INT]{Introduction.}

\noindent
In this paper we study the equicontinuity properties of meromorphic mappings with
values in compact complex surfaces.
\smallskip
Let $\dd_r$ be disk of radius 
$r$ in $\cc$, $\dd:=\dd_1,$ $\dd^2_r$ the
bidisk of radius $r$ of center $0$ in $\cc^2.$ Given $0 < r_1 < r_2,$
$R^2(r_1,r_2)$ denotes the
crown $\dd^2_{r_2} \backslash \overline{\dd}^2_{r_1}$. Recall further that
the Hartogs figure  $H^2_{\epsilon}$ in $\cc^2$ is the following domain:

\begin{equation}
\eqqno(hart)
H^2_{\epsilon}= \left[ \dd_{\epsilon}\times \dd\right] \cup\left[
\dd\times(\dd\backslash \overline\dd_{1-\epsilon} ) \right],
\end{equation}
where  $0< \epsilon <1$. 

The main result of this paper is the following.

\begin{thm}Let ${\cal F}$ be a family of meromophic mappings from $\dd^2$ to a
compact complex surface $X$. Suppose that there exists $ \epsilon > 0$  such that 
every  $f \in {\cal F}$ is holomorphic on the Hartogs figure $H^2_{\epsilon}\subset \dd^2$ 
and $\calf $ is equicontinuous  on $H^2_{\epsilon}$. 
Then for every  $ 0< r <1,$ there is a constant $C_r$ such that 
\begin{equation}
\vol\left(\Gamma_{f}\cap(\dd^2_{r}\times X)\right)\le C_r,
\end{equation}
for every $f \in {\cal F}.$

\end{thm}

\smallskip 
Here $\Gamma_{f}$ denotes the graph of $f$. Recall that a meromorphic mapping $f$
between complex manifolds $Y$ and
$X$ is defined by its graph $\Gamma_f,$  which is an analytic subset of the product
$Y\times X$, satisfying
the following conditions:

\smallskip
\textit{1.} $\Gamma_f$ is a locally irreducible analytic subset of $Y\times X$;

\smallskip
\textit{2.} The restriction $\pi\arrowvert \Gamma_f : \Gamma_f \longrightarrow Y $ of
the natural projection $\pi: Y \times X \longrightarrow Y$ to $\Gamma_f$ is
proper, surjective and generically one to one. A family ${\cal F}$ of mappings
$Y\longrightarrow Z,$ where $Y$ and $Z$ are  metric spaces, is called equicontinuous
if for every $\epsilon >0,$ there is
$\delta>0$ such that $dist\left(f(z),f(w) \right)< \epsilon $
for every $f\in {\cal F}$ and  $dist(w,z)<\delta.$   

\smallskip 
One of the main ingredients of the proof of Theorem 1 is the following result.

\begin{prop} Let ${\cal F}$ be a family of meromophic mappings from $\dd^2$ to a
compact complex surface $X$. Suppose that there exists $ \epsilon > 0$  such that 
every  $f \in {\cal F}$ is holomorphic in $R^2(1-\epsilon,1)$ and the family ${\cal
F}$ is equicontinuous on $R^2(1-\epsilon ,1)$. Then  there exists
a constant $C_r$ such that
\begin{equation}
\vol\left(\Gamma_{f}\cap(\dd^2_{r}\times X)\right)\le C_r,
\end{equation}
for every $f \in {\cal F}.$

\end{prop}

In \cite{Iv1} the notions of weak and strong convergence of sequences of meromorphic
mappings were introduced. We recall them in section 4 and prove the following.
\begin{prop}
If $X$ is a non-projective compact complex surface, then every weakly convergent
sequence of meromorphic 
mappings from a bidisk to $X$ converges strongly.
\end{prop}

\begin{rema} {\bf 1.} \rm
Remark that neither  Theorem 1 nor Corollary 1 follow from the Oka-type estimates
proved in \cite {F-S}. The reason is that to estimate the volume of $\Gamma_{f}$ one
needs, in particular, to estimate the integral
$\int_{X}\left(f^\ast\omega\right)^2,$ and the current $\left(f^\ast\omega\right)^2$
is of bidimension $(0,0).$

\smallskip\noindent{\bf 2.}
  In the case $X$ is a K\"ahler manifold, Theorem 1  is proved In \cite{Iv1}.

\smallskip\noindent{\bf 3.}
In Section 5, we show by an counterexample that Theorem 1  is not valid in general  when $dimX
\geq 3.$
\end{rema}

\smallskip
I'm grateful to my advisor S. Ivashkovich for guiding me in this research work.

\newsect[sect.AREA]{Estimates of areas of sections}
 Let $h$ be a hermitian metric on a complex manifold $X,$ and let $\omega$ be the
$(1,1)$-form canonically associated with $h$.
Metric $h$ (and form $\omega$) is called pluriclosed (or $dd^c$-closed) if
$dd^c\omega=0$. By \cite{Ga}, every compact complex surface admits a pluriclosed
metric form. For $f: \dd^2 \longrightarrow X $ a meromorphic map, we denote by $A_
{f} \subset  \dd^2(r)$ the  set of points of inderterminacy  of $f$.
Consider the  current $T_{f}=f^\ast\omega$ on $\dd^2$.
Write $$ {T_f}=\frac{i}{2 } t^{\alpha\bar{\beta}}_f dz_\alpha \wedge d\bar{z_\beta},
$$ where $ t^{\alpha\bar{\beta}}_f$ are distributions on $\dd^2,$ smooth on
$\dd^2\backslash A_{f}.$

A complex manifold $X$ is called disc-convex if for every compact set $K\comp X$
there exists another
compact set $\hat K$ such that for every analytic disc $\phi:\overline \dd\to X,$
with $\phi (\d\dd)\subset K,$ 
one has $\phi(\overline\dd)\subset \hat K$. Note that compact manifolds and Stein
manifolds are disc-convex. More generally,  each 1-convex manifold is disc-convex.

\begin{prop}
Let ${\cal F}$ be a family of meromophic mappings from $\dd^2$ to a disc-convex
manifold $X,$ which admits a pluriclosed metric form. Suppose that for some  $0 <
\epsilon < 1,$ the family ${\cal F}$ is holomorphic and equicontinuous on $R^2(\epsilon ,1).$  Then
for every $0 < r < 1,$ areas of graphs of restrictions  $\Gamma_f \cap \left(
\dd_{z_1}(r) \times X \right) $ of $f$ to the discs $\dd_{z_1}(r)= \left\lbrace
z_1\right\rbrace \times \dd(r)$ are uniformly bounded in $z_1\in\dd(r)$ and
$f\in\calf$.
\end{prop}
\proof The proof will be done in three steps. First two we shall state in the
form of a lemma.

\begin{lem}
Distributions $t^{\alpha\bar{\beta}}_f $ are locally integrable in $\dd^2.$
\end{lem}

\proof\ \rm
Note that the family of smooth forms $\{T_{f}\lvert_{R^2(r,1)}:f \in{\cal F}\}$ is
equicontinuous on $ R^2(r,1).$
Fix $r_1 \in [ r,1 [$ close enough to 1.  Let $z=(z_1,z_2)\in \dd^2.$ 

\noindent
Set $\dd_{z_1}\left( r_1 \right)= \left\lbrace
z_1\right\rbrace \times\dd\left(r_1\right).$
 Consider functions $a_{f}$ given by 
\begin{equation}
a_{f}(z_1)= \int_{\dd_{z_1}(r_1)}T_{f} =\frac{i}{2 }\int_{\dd_{z_1}(r_1)}
t^{2\bar{2}}_{f} d{z_2} \wedge d\bar{z_2}.
\end{equation}
Functions $a_{f}$ are well-defined and smooth on
${\dd\backslash\pi(A_{f})},$ where $\pi :
\dd^2\longrightarrow \dd$ is the canonical projection on the first factor. The
condition that $dd^cT_f=0$ implies, in particular, that
\begin{equation}
\frac {\partial^2 t^{2\bar{2}}_{f}}{\partial{z_1}\partial {\bar{z_1}}}+ \frac{
\partial^2 t^{1\bar{1}}_{f}}{\partial{z_2}\partial{\bar{z_2}}}-\frac{ \partial^2
t^{1\bar{2}}_{f}}{\partial{z_2}\partial{\bar{z_1}}}-\frac{\partial^2
t^{2\bar{1}}_{f}}{\partial{z_1}\partial{\bar{z_2}}}=0.
\end{equation}
 Now we can estimate the Laplacian of $a_{f}$ on ${\dd \backslash\pi(A_{f})}.$ We have 
\[
\Delta a_{f}(z_1)= \frac{i}{2} \int\limits_{\dd_{z_1}(r_1)} \frac {\partial^2
t^{2\bar{2}}_{f}}{\partial{z_1}\partial {\bar{z_1}}} d{z_2} \wedge d\bar{z_2}=
\]
\[
= \frac{i}{2}\left( \int\limits_{\dd_{z_1}(r_1)}\frac{ \partial^2
t^{1\bar{2}}_{f}}{\partial{z_2}\partial{\bar{z_1}}}+ \frac{ \partial^2
t^{2\bar{1}}_{f}}{\partial{z_1}\partial{\bar{z_2}}}- \frac{ \partial^2
t^{1\bar{1}}_{f}}{\partial{z_2}\partial{\bar{z_2}}}\right) d{z_2} \wedge d\bar{z_2}=
\frac{i}{2} \int\limits_{\partial\dd_{z_1}(r_1)} \frac{ \partial
t^{1\bar{2}}_{f}}{\partial{\bar{z_1}}} d\bar{z_2}
\]

\begin{equation}
 + \frac{i}{2} \int\limits_{\partial\dd_{z_1}(r_1)} \frac{ \partial
t^{2\bar{1}}_{f}}{\partial{{z_1}}} d{z_2}
-\frac{i}{2} \int\limits_{\partial\dd_{z_1}(r_1)} \frac{ \partial
t^{1\bar{1}}_{f}}{\partial{\bar{z_2}}} d\bar{z_2}=:\varphi_{f}\left( z_1\right).
\end{equation}

Since $T_{f}$ is smooth in a neighborhood of $\dd\times \partial \dd,$ 
$\varphi_{f}$ are smooth in the whole unit disc $\dd,$ for every $f \in  {\cal F}.$ 
Set $\psi_{f}(z_1)=\varphi_f* ln\lvert \zeta-z_1\lvert.$ Then
\begin{equation}
\Delta \psi_{f}=\varphi_{f}.
\end{equation}

\noindent
Set
\begin{equation} 
h_{f}:= a_{f}-\psi_{f}.
\end{equation}
 Since $a_{f}$ is positive on ${\dd \backslash\pi(A_{f})}$ and $\psi_{f}$ is smooth on
$\dd$, $h_{f}$ is bounded from below on
$\dd.$ Also $\Delta h_{f}=0$ on
${\dd \backslash\pi(A_{f})}$. Therefore $h_{f}$ extends to a
superharmonic functions on $\dd.$ Therefore $h_{f}\in
L^1_{loc}(\dd),$ see \cite[2.5 Theorem 1]{R}. It follows that $ t^{2\bar{2}}_{f} \in
L^1_{loc}\left(\dd^2\right).$ A similar argument shows that $t^{1\bar{1}}_{f}$ is
locally integrable in $\dd^2.$ Positivity of $T_{f}$ implies that $
t^{1\bar{1}}_{f} t^{2\bar{2}}_{f}\geqslant\vert
t^{1\bar{2}}_{f}\vert^2.$ So, in particular, 

\begin{equation}
\int\limits_{\dd^2(r_1)} \vert
t^{1\bar{2}}_{f}\vert\leqslant\int\limits_{\dd^2(r_1)}\sqrt{t^{2\bar{2}}_{f}}.
\sqrt{t^{1\bar {1}}_{f}}\leq \sqrt{ \int\limits_{\dd^2(r_1)}t^{2\bar{2}}_{f}}.\sqrt{
\int\limits_{\dd^2(r_1)}t^{1\bar{1}}_{f}},
\end{equation}
 which gives that $t^{\alpha\bar{\beta}}_f \in L^1_{loc}\left( \dd^2\right).$

\smallskip\qed

\begin{lem}
Under the hypotheses of Proposition 3, the family $\left\lbrace t^{\alpha\bar{\beta}}_f
\right\rbrace_{f\in \cal F} $ is  uniformly bounded in $L^1\left( \dd^2(r)\right),$
for every $0 < r <1.$ 
\end{lem}

\proof\ \rm
Let $\dd^2(r)\subset \dd^2(r_1)\subset \dd^2.$ According to (2.3),  we have
$$\Delta a_{f}(z_1)=\varphi_{f}\left( z_1\right),$$ where $\varphi_{f}$ is a smooth
function on $\dd.$ Moreover from (2.3) we see that the family $\varphi_f$ is
equicontinuous.  Fonction $h_f$ given by  $(2.5)$ is superharmonic in $\dd$ and
harmonic on ${\dd \backslash\pi(A_{f})}.$ Therefore
\begin{equation}
\Delta h_f=-\sum_{z_j \in \pi(A_f)} c_j(f)\delta_{z_j(f)} = :\mu_f,
\end{equation}
where $c_j>0.$  Therefore, we can rewrite  (2.5)  as 
\begin{equation}
\Delta a_f= \Delta \psi_f -\sum_{z_j \in \pi(A_f)} c_j(f)\delta_{z_j(f)}.
\end{equation}
Note that $\left\lbrace \psi_f\right\rbrace_{f\in \cal F} $ are equicontinuous on
$\dd(r).$

Fix $\delta > 0$ such that  $\dd(\delta,z_j),$ which $z_j \in \pi(A_f),$ are
pairwise disjoint. Let $\phi$ be a test function on $\dd(r)$ with $supp\phi
\subset\subset \dd(\delta,z_j).$ 

\noindent 
Let $a^{\epsilon}_f(z_1)=\frac{i}{2 } \int\limits_{\dd_{z_1}(r_1)}
t^{2\bar{2}}_{f,\epsilon} d{z_2} \wedge d\bar{z_2},$ where
$t^{2\bar{2}}_{f,\epsilon}$ is the smoothing of $t^{2\bar{2}}_{f}$ by convolution.
Since $t^{2\bar{2}}_{f,\epsilon} \longrightarrow t^{2\bar{2}}_{f}$ in  $L^1\left(
\dd^2(r)\right),$ we get by the Fubini Theorem  that $ a^{\epsilon}_f
\longrightarrow a_f$ in $L^1\left(\dd(r)\right).$ Hence, 
\[
< \Delta a^{\epsilon}_f,\phi > = \int\limits_{\dd(\delta, z_j)} a^{\epsilon}_f(z_1)
\Delta \phi(z_1) dz_1\wedge d\bar z_1 =\frac{i}{2 }\int\limits_{\dd(\delta,
z_j)}\int\limits_{\dd_{z_1}(r_1)} t^{2\bar{2}}_{f,\epsilon}\Delta \phi(z_1) d{z_2}
\wedge d\bar{z_2}\wedge dz_1\wedge d\bar z_1 =  
\]

\[
= \frac{i}{2 } \int\limits_{\dd_{z_1}(r_1)} \int\limits_{\dd(\delta, z_j)}
t^{2\bar{2}}_{f,\epsilon}\Delta \phi(z_1) dz_1\wedge d\bar z_1 \wedge d{z_2} \wedge
d\bar{z_2}=
\]
\[
= \frac{i}{2 } \int\limits_{\dd_{z_1}(r_1)} \int\limits_{\dd(\delta, z_j)}
\frac{\partial^2 t^{2\bar{2}}_{f,\epsilon}}{\partial z_1 \partial \bar z_1}
\phi(z_1)dz_1\wedge d\bar z_1 \wedge d{z_2} \wedge d\bar{z_2}= 
\]
\[
=\frac{i}{2 }\int\limits_{\dd(\delta, z_j)} \phi (z_1) \left(
\int\limits_{\dd_{z_1}(r_1)} \frac{\partial^2 t^{2\bar{2}}_{f,\epsilon}}{\partial
z_1 \partial \bar z_1} d{z_2} \wedge d\bar{z_2} \right) dz_1\wedge d\bar z_1. 
\] 
Using (2.5), we obtain
\[
< \Delta a^{\epsilon}_f,\phi >= \frac{i}{2 }\int\limits_{\dd(\delta, z_j)} \phi
(z_1) \left( \frac{i}{2} \int\limits_{\partial\Delta_{z_1}(r_1)} \frac{ \partial
t^{1\bar{2}}_{f,\epsilon}}{\partial{\bar{z_1}}} d\bar{z_2} + \frac{i}{2}
\int\limits_{\partial\dd_{z_1}(r_1)} \frac { \partial
t^{2\bar{1}}_{f,\epsilon}}{\partial{{z_1}}} d{z_2} \right)dz_1\wedge d\bar z_1
\]

\[
- \frac{i}{2 }\int\limits_{\dd(\delta, z_j)} \phi (z_1) \left( \int\limits_{\partial
\dd_{z_1}(r_1)} \frac{ \partial t^{1\bar{1}}_{f,\epsilon}}{\partial{\bar{z_2}}}
d\bar{z_2}\right) dz_1\wedge d\bar z_1 \longrightarrow < \varphi_{f},\phi >,
\]
as $\epsilon \longrightarrow 0.$

Therefore, $\Delta a_f= \varphi_f $ in $\dd(r)$ in the sens of distributions. By
regularity of the Laplacian, $a_f \in \cal C^{\infty}.$ 
It follows that $ t^{2\bar{2}}_{f}$ is uniformly bounded in $L^1\left(
\dd^2(r)\right).$ Same is true for $
t^{1\bar{1}}_{f}.$ Using (2.6), we see that  $t^{\alpha\bar{\beta}}_f$ are 
uniformly bounded  in $L^1_{loc}\left( \dd^2(r) \right).$

\smallskip\qed        
\begin{rema}\rm
This lemma can be proved also using the Oka-type inequality from \cite{F-S}.
\end{rema}

\smallskip\noindent{\slsf End of the proof.}
We  denote by  $\hat{f}$ the mapping $\hat{f}(z)=(z,f(z))$ into the graph  $\Gamma_f$.
Fix some $0<r<1$. 
The area of  $ \Gamma_f \cap \left( \dd_{z_1}(r) \times X \right) $
is given by
\[
area \left( \Gamma_{f} \cap \left( \dd_{z_1}(r) \times X \right) \right)= area\hat{f}(
\dd_{z_1}(r))=  
\int\limits_{\dd_{z_1}(r)} \left( T_{f}\arrowvert \dd_{z_1}(r)+ dd^c\lvert
z_2\lvert^2\right) =a_f \left(z_1 \right) + \int\limits_{\dd_{z_1}(r)} dd^c\lvert
z_2\lvert^2.
\]
 According to Lemma 2,  $a_f$ are uniformly bounded. Therefore, the areas of $
\Gamma_f \cap \left( \dd_{z_1}(r) \times X \right) $ are uniformly bounded in $f$ and $z_1$. The proof of Proposition
3 is complete.

\smallskip\qed

\newsect[sect.VOL]{Estimates of volumes }

In this section we shall prove the Proposition 1 stated in Introduction. Along the
proof,  we shall crucially 
use the following  results of Barlet (see \cite{Ba1} and \cite{Ba2}). We recall that
a striclty q-convex function $\rho$ on the complex space $X$ with $dimX=N$ is a real
valued ${\cal C}^2$-function such that the hermitian matrix consisting of the
coefficients of the $(1,1)$-form $dd^c\rho$ has at least $N-q+1$ positive
eigenvalues at all points of $X.$  The complex space $X$ is called $q$-complete if
it admits a strictly $q$-convex exhaustion function $\rho$: $X\longrightarrow
\rr^+$, $q\ge 1$. Note that, in the case $q=1,$ $X$ is Stein.
\smallskip

\smallskip
\medskip\noindent{{\bf B1. }\it Let $C$ be a $q$-dimensional compact analytic
subspace of a complex space $X$. Then 
$C$ admits a fundamental system of $(q+1)$-complete neighborhoods.}

\medskip\noindent{{\bf B2. }\it Let $X$ be a $(q+1)$-complete complex space and let
$\rho : X\to \rr^+$ be a strictly 
$q$-convex function. Let $h$ be some ${\cal C}^2$-smooth hermitian metric on $X$.
Then there exists a hermitian 
metric $h_1$ on $X$ and a function $c:\rr^+\to \rr^+$ (both of class ${\cal C}^2$)
such that

\noindent

\textit{(i)} $h_1\geq h$;
\smallskip

\textit{(ii)} the $(q+1,q+1)$-form $\Omega = dd^c[(c\circ \rho)\omega_{h_1}^{q}]$ is
strictly positive on $X$.}
\smallskip

\smallskip
Here $\omega_h$ is the $(1,1)$-form canonically associated with $h$. A
$(q+1,q+1)$-form $\Omega $ is called strictly positive if for any $x \in X$ and
linearly independent vectors $v_1,...,v_{q+1} \in T_{x}X$ one has $ \Omega_x\left(
iv_1\wedge \bar v_1,..., iv_{q+1}\wedge \bar v_{q+1} \right)> 0.$

It should be noted that  Enriques-Kodaira classification of compact complex 
surfaces implies that if a compact surface $X$ contains an infinite number of rational
curves, then $X$ is projective and, in particular,  K\"ahler, (see \cite{BHPV}).
Theorem 1  for K\"ahler manifold  was proved in  \cite{Iv1}. Therefore, in our proof
we can suppose that our surface $X$ contains at most finite number of rational
curves.

\smallskip 
We need to prove the uniform estimate of volumes of $\Gamma_{f}$
\begin{equation}
 \vol\left( \Gamma_{f} \right)=\int\limits_{\dd^2} \left( T_{f}+ dd^c\Vert z
\Vert^2\right)^2,  (f \in \cal F).
\end{equation}
If the result is false, then there exists a sequence ${f_n \in {\cal F}}$ of
meromorphic mappings from $\dd^2$ to a compact complex surface  $X$ holomorphic on
$R^2(\epsilon,1)$ such that
\smallskip

\smallskip
\textit{1} $f_n \lvert_{R^2(\epsilon,1)}$ uniformly converge  to $f$;
\smallskip

\textit{2} $\vol\left(\Gamma_{f_n} \right)\longrightarrow \infty.$
\smallskip

\textit{3} Taking a subsequence we can suppose that $\Gamma_{f_n} $ converge in
Hausdorff metric on compacts subsets in $\dd^2 \times X$
 to a closed set $\Gamma$.
\smallskip

\smallskip 
Since $f_n(\partial\dd^2)$ are homologous to zero in $X$ for every n, 
$f(\partial\dd^2)$ is homologous to zero in $X$ . Therefore $f$ extends
meromorphically to $\dd^2$, see Lemma 2.5 in \cite {Iv2}. So we have  $\Gamma\cap \left[\left( \dd^2\backslash
A_{f} \right) \times X \right]=\Gamma_{f}.$ We  prove that $\Gamma_{f_n}\cap
\left[\left( \dd^2\backslash A_{f} \right) \times X \right]\longrightarrow
\Gamma_{f}\cap \left[\left( \dd^2\backslash A_{f}  \right) \times X \right]$.
Indeed, if $z_0\notin A_{f}$, i.e., $f$ is holomorphic in a neighborhood of $z_0,$
then we can find  neighborhoods $U \ni z_0$ and  $V \ni y_0=f(z_0)$ such that
$f(\bar U)\cap \partial V= \emptyset$. But then, for enough large $n,$ $f_n(\bar
U)\cap \partial V= \emptyset.$ Let $i: V\longrightarrow V'$ be an isomorphism onto
the bounded open set in $\cc^N$. Now  $i\circ f_n: U\longrightarrow V'$ is a
sequence of holomorphic mappings of $U$ into $V'$. So we can find a convergent
subsequence. Therefore, the graph of the limit must coincide with $\Gamma_f \cap (U\times V).$

\smallskip
We can suppose that $A_{f}=\left\lbrace 0 \right\rbrace $. If $(0,x) \in \left[
\Gamma \cap \left( \left\lbrace 0 \right\rbrace \times X \right) \right]\backslash
\Gamma_{f}$ then there exists a rational curve $C_x$ containing $x$. Indeed, if 
$(z_n,x_n) \in \Gamma_{f_n}$ converge to $(0,x),$ then by  to Proposition 3,
$area(f_n(\Delta_{z_n}))$ is uniformly bounded. By Bishop's Theorem (see for example
\cite{St}), we can assume, after passing to a subsequence, that $f_n(\Delta_{z_n}) $
 converge to $f(\dd_{0})\cup C,$ where $C$ is a finite union of rational curves, see Lemma 10 in \cite{Iv3}. But $(0,x) \notin f(\dd_{0})$ because  $(0,x) \notin
\Gamma_f.$  The only possibility is that $(0,x) \in C$. Since the number of rational
curves in $X$ is finite, we see that 
$$\Gamma= \Gamma_{f}\cup\left( \left\lbrace \ 0 \right\rbrace \times \cup_{ i=1}^{N}
C_i \right).$$

According to \cite{Ba1},  $\Gamma \cap \left( \left\lbrace 0 \right\rbrace \times X
\right)$ admits a 2-complete neighborhood $W\subset  \left\lbrace 0 \right\rbrace
\times X.$ In addition  $\dd^2 \times W$ is also  2-complete. We apply Barlet
theorem \cite{Ba2}, by taking $\rho$ to be  strictly 2-convex exhaustion of $\dd^2
\times W$ in order to have a strictly $dd^c$-exact (2,2)-form $\Omega$ on $\dd^2
\times W$. Let $\tau $ be a fixed (1,1)-form of class ${\cal C}^2$ such that
$dd^c\tau=\Omega$. We can suppose that
$\tau \in {\cal C}^2( \overline \dd^2 \times W),$ i.e, $\tau$ is smooth up to the
boundary. Hence,

\begin{equation}
\vol\left(\Gamma_{f_n}\right)\lesssim \int\limits_{\Gamma_{f_n}\arrowvert \left(
\dd^2 \times W \right)} \Omega = \int\limits_{\Gamma_{f_n}\arrowvert \left(   \dd^2
\times W \right)} dd^c \tau=
\end{equation}
\begin{equation}
= \int\limits_{\Gamma_{f_n}\arrowvert \partial \left( \dd^2 \times W \right)} d^c
\tau\leq C,
\end{equation}
where the constant $C$ does not depend on $f_n.$  Because $f_n$ converge on compacts
subsets outside of zero,  $d^c\tau$ is of class ${\cal C}^1$ on $\overline \dd^2
\times W,$ which is a contradiction.  The proof of Proposition 1  is complete.

\newsect[sect.COROL]{Proof of Theorem 1 } 

Now  we prove Theorem 1. If we proceed by contradiction. Then, there exists a
sequence ${f_n \in {\cal F}}$ of meromorphic mappings  from $\dd^2$ to a compact
complex surface $X$ holomorphic on $H^2(\epsilon)$ such that
\smallskip

\smallskip
\textit{1} $f_n \lvert_{H^2(\epsilon)}\rightrightarrows f;$
\smallskip

\textit{2} $\vol\left(\Gamma_{f_n}\arrowvert_{  \dd \times \Delta }
\right) \longrightarrow \infty$.
\smallskip

\smallskip
According to \cite[Theorem 1]{Iv2}, $f$ extends to $\dd^2\backslash A,$ where
$A$ is discrete. Take $s_0\in A$ and let $S^3_{s_0}(r)$ be some euclidean sphere
centered at $s_0$ such that $S^3_{s_0}(r)\cap A=\emptyset.$ Since  $f_n(S^3_{s_0})$
is homologous to zero in $X,$ for $n\geq 1,$  $f(S^3_{s_0})$ is
homologous to zero in $X$. According to \cite{Iv2}, this implies that $f$ extends
onto the ball $B_{s_0}$ with $\partial B_{s_0}= S^3_{s_0}.$ Therefore, by
\cite[Proposition 1.1.1]{Iv1}, $f_n\longrightarrow f$ on  compacts subsets of
$\dd^2\backslash A_{f}.$ An application of Proposition 1 applied to $f$ gives a
contradiction.

\smallskip\qed

\newsect[sect.COROL]{Proof of Proposition 2}

Let $\left\lbrace f_n \right\rbrace$ be a sequence of meromorphic mappings from $\dd^2$ to a complex compact surface $X.$ Let us recall the following definitions from \cite{Iv1}.
	
\begin{defi}
We  say that $\left\lbrace f_n \right\rbrace $ converges strongly (s-converge) on
compacts subsets  in $\dd^2$ to a meromorphic map $f: \dd^2\longrightarrow X$ if for
any compact set $K \subset X$ 
$$ {\cal H}- \lim\limits_{n \rightarrow \infty} \Gamma_{f_n}\cap \left( K \times X
\right) = \Gamma_{f}\cap \left( K \times X \right).$$ 

Here ${\cal H}- \lim\limits$ denotes the limit in the Hausdorff metric.
\end{defi}

\begin{defi}
We say that a sequence of meromorphic mappings $f_n: \dd^2\longrightarrow X $
converges weakly (w-converge) on $\dd^2$ to a meromorphic map $f:\dd^2
\longrightarrow X$ if there exists a discrete  subset $A$ of $\dd^2$ such that  $f_n
\longrightarrow f$ strongly on compacts subsets in  $\dd^2\backslash A.$
\end{defi}

\smallskip\noindent{\slsf Proof of Proposition 2. } Let $f_n: \dd^2\longrightarrow X $ be a sequence
of meromorphic mappings $f_n$ from  $\dd^2$ to compact surface $X$ converging weakly
in $\dd^2$ to a meromorphic map $f: \dd^2 \longrightarrow X$. Therefore, there
exists a discrete subset $A$ of $\dd^2$ such that $f_n\longrightarrow f$ in
$\dd^2\backslash A.$ Without loss of generality we can suppose that $A=\{0\}$. 
We shall prove that, for every compact $K\subset \dd^2$
$$ {\cal H}- \lim\limits_{n \rightarrow \infty} \Gamma_{f_n}\cap \left( K \times X
\right) = \Gamma_{f}\cap \left( K \times X \right),$$ i.e., the sequence of graphs $
\Gamma_{f_n}\cap \left( K \times X \right)$ converges in the Hausdorff metric to the
graph of the limit. We have earlier proved that $\vol \left(  \Gamma_{f_n}\right) $
are bounded. Therefore,  there exists a subsequence, (still denoted by
$\Gamma_{f_n}\cap \left( K \times X \right)$), which converges to 

\[
\Gamma = \Gamma_{f} \cap \left( K \times X \right) \cup \left( \left\lbrace 0
\right\rbrace \times X\right).
\]
We have already proved in the proof of Proposition 1 that if $x\in \Gamma \backslash
\Gamma_f,$ then there exists a rational curve $ C\ni x.$ Since $X$ is not  K\"ahler
compact surface it contains at most a finite number of rational curves, which is a
contradiction. Therefore $\Gamma=\Gamma_f$.

\smallskip\qed

\newsect[sect.EXMP]{Example}

\noindent
In general, if the manifold $X$ is of dimension bigger than two, Theorem 1 doesn't
hold true. Consider  metric form on a Hopf manifold $  X^n=\left(
\cc^n\backslash\left\lbrace 0\right\rbrace \right) /\left(z\sim 2z \right)$ of any
dimension $n\geq 2$ defined by
\begin{equation}
\omega=\frac{i}{2}\frac{(dz,dz)}{\parallel
z\parallel^2}=\frac{i}{2}\frac{dz_1\wedge d\bar z_1 +...+ dz_n\wedge
d\bar z_n  }{\parallel z\parallel^2}.
\end{equation}

We use the following notations: $(dz,dz)= dz_1\wedge d\bar z_1 +...+ dz_n\wedge
d\bar z_n$, $(dz,z)= dz_1\wedge \bar z_1 +...+ dz_n\wedge
\bar z_n$ and $(z,dz)= z_1\wedge d\bar z_1 +...+ z_n\wedge d\bar z_n.$
\noindent

\begin{exmp}\rm
\noindent
Consider holomorphic maps $f_n:B^2\longrightarrow X^3 $  defined by
$f_n(z _1,z_2)= \left( z_1,z_2,\frac{1}{n}\right).$ So, $\left\lbrace
f_n\right\rbrace $ is equicontinuous on $R^2(\epsilon,1),$ for any $\epsilon > 0.$
But as $n \longrightarrow \infty,$  $\vol\left(\Gamma_{f_n}\right) \longrightarrow
\infty.$
Indeed, if $z'=(z_1,z_2 ),$ then

$$\vol\left(\Gamma_{f_n}\right) \geq \int\limits_{f_n(B^2)}\omega^2
=\int\limits_{B^2}(f_n^*\omega)^2 =$$

$$=\int\limits_{B^2(1)} \left(\frac{i}{2}\frac{ dz_1\wedge d\bar z_1+ dz_2\wedge
d\bar z_2 }{\lvert z_1\lvert^2+\lvert z_2\lvert^2+ \frac{1}{n^2}}\right)^2
=2 \int\limits_{B^2(1)} \left( \frac{i}{2}\right)^2 \frac{ dz_1\wedge d\bar
z_1\wedge dz_2\wedge
d\bar z_2 }{\lvert z_1\lvert^2+\lvert z_2\lvert^2+ \frac{1}{n^2}} \approx$$
$$\approx \int\limits_{0}^{1} \frac{r^3 dr } {\left( r^2 +
\frac{1}{n^2}\right)^2}\longrightarrow +\infty,
$$
as $n\longrightarrow \infty.$
\end{exmp}

\smallskip\qed

\begin{defi}
We say that $h$ is plurinegative if $dd^c\omega \leq 0$.
\end{defi}

\smallskip Remark that 
if $n=2,$ then for $\omega$ defined by (6.1) one has  $dd^c\omega=0,$ i.e, $\omega$ is pluriclosed.

\begin{lem} If $n\geq 3,$ then $\omega$ is plurinegative but not pluriclosed.

\end{lem}
\proof \rm Let $\omega$ is the metric form on $X$ defined as in (6.1). Then
\noindent
\[
\bar\partial\omega= -\frac{i}{2}\frac{(z,dz)\wedge(dz,dz)}{\parallel
z\parallel^4 }
\]
\noindent
\[
\partial\bar\partial\omega= -\frac{i}{2}\frac{(dz,dz)\wedge(dz,dz)}{\parallel
z\parallel^4}+i\frac{(dz,z)\wedge(z,dz)\wedge(dz,dz)}{\parallel
z\parallel^6}
\]
\noindent
\[
dd^c\omega=2i\partial\bar\partial \omega=\frac{(dz,dz)}{\parallel
z\parallel^6}\big(\parallel z \parallel^2(dz,dz)-2(dz,z)\wedge(z,dz)\big).
\]
If $v=(v_1,...,v_n) \in \cc^n$, then
\[
 (\parallel z \parallel^2(dz,dz)-2(dz,z)\wedge(z,dz))\wedge
 (v,\bar v)=\parallel z \parallel^2\parallel
 v\parallel^2-2(v,z)(z,v)
 \]
 \[
 =\parallel z \parallel^2\parallel
 v\parallel^2-2|(v,z)|^2\geq 0,
 \]
\noindent
according to the Schwarz inequality. Therefore,

$$dd^c\omega(ib_1,\bar b_1,ib_2,\bar b_2)=-\parallel b_1 \parallel^2(\parallel z
\parallel^2\parallel b_2 \parallel^2-2\lvert(b_2,z)\lvert^2) \leq 0,$$
i.e., $dd^c\omega$ is non positive and not identically zero.
\noindent
 Hence the result.
\smallskip\qed
\begin{rema}\rm
Note that if $n \geq 3,$ $X^n$ admits no pluriclosed metric form. It is sufficent to
prove this for $n=3.$ If $\Omega$ is such a form, then for $\omega $ as in (6.1)
one has 
$$ 0 > \int \Omega\wedge dd^c\omega = \int dd^c\Omega \wedge \omega =0,$$ which is
clearly impossible.
\end{rema}
\begin{rema}\rm
We believe that our Theorem holds for $X$ of any dimension admitting a pluriclosed
metric form. The result should be true also for Hartogs figures in all dimensions. 
\end{rema}

\def\entry#1#2#3#4\par{\bibitem[#1]{#1}
{\textsc{#2 }}{\sl{#3} }#4\par\vskip2pt}

\end{document}